\documentclass[a4paper,10pt]{article}
\usepackage[utf8]{inputenc}
\usepackage{amssymb,latexsym}
\usepackage{amsmath,amsthm}
\usepackage{graphicx}
\usepackage{color,epsfig}
\usepackage{multicol}
\usepackage{multirow}

\newtheorem{theorem}{Theorem}
\newtheorem{lemma}{Lemma}

\def\qed{\ifhmode\unskip\nobreak\fi\quad\ifmmode\Box\else$\Box$\fi}

\long\def\comment#1{}

\def\cro{{\mbox {\sc cr}}}
\def\ocr{{\mbox {\sc ocr}}}
\def\pcr{{\mbox {\sc pcr}}}

\title{Crossing lemma for the odd-crossing number}
\author{J\'anos Karl\\
\small  Budapest University of Technology and Economics \\[-0.8ex]
\small \texttt{karlj@math.bme.hu}\\
and\\
G\'eza T\'oth\thanks{Supported by National Research, Development and Innovation Office, NKFIH,
K-131529 and ERC Advanced Grant "GeoScape" 882971.}\\
\small Alfr\'ed R\'enyi Institute of Mathematics, Budapest
\small and
\small Budapest University of Technology and Economics \\[-0.8ex]
\small \texttt{geza@renyi.hu}}

\begin{document}

\maketitle

\begin{abstract}
  A graph is $1$-planar, if it can be drawn
  in the plane such that there is at most one crossing on 
every edge.
It is known, that $1$-planar graphs have at most $4n-8$ edges.

We prove the following odd-even generalization. If a graph can be drawn 
in the plane such that
every edge is crossed by at most one other edge {\em an odd number of times},
then it is called 1-odd-planar and it has at most $5n-9$ 
edges. As a consequence, we improve the
constant in the Crossing Lemma for the odd-crossing number,
if adjacent edges cross an even
number of times.
We also give upper bound for the number of edges of $k$-odd-planar graphs.

\end{abstract}

\section{Introduction}

By a {\em graph} we always mean a {\em simple graph}, that is, a graph with no loops and no parallel edges.
We use the term {\em multigraph} if loops and parallel edges are allowed.
A {\em drawing} of a (multi)graph in the plane is a representation
such that vertices are represented by  
distinct points and its edges by curves connecting 
the corresponding points. We 
assume that no edge passes through any vertex
other than its endpoints, no two edges touch each other
(i.e., if two edges have a common interior point, then at this
point they properly cross each other), no three edges
cross at the same point, and two edges cross only finitely many times. 

The {\em crossing number} of a graph $G$, $\cro(G)$, is the minimum number of crossings 
(crossing points)
over all drawings of $G$. 
The {\em pair-crossing number}, $\pcr(G)$, is the minimum number of pairs of crossing edges
over all drawings of $G$.  
In an optimal drawing for $\cro(G)$, any two edges cross at most once \cite{S04}.
Therefore, 
it is not easy to see the difference between these two definitions.
Indeed, there was some confusion in the literature between these two notions, 
until the systematic study of their relationship \cite{PT00a}. 
Clearly, $\pcr(G)\le\cro(G)$, and 
in fact, we cannot rule out the possibility, that $\cro(G)=\pcr(G)$ 
for every graph $G$. Probably 
it is the most interesting open problem in this area.
From the other direction, the best known bound is 
$\cro(G)=O(\pcr(G)^{3/2}\log\pcr(G))$ \cite{S17}, \cite{KT21}.




The {\em odd-crossing number}, $\ocr(G)$, is the minimum number 
of pairs of edges that cross an odd number of times,
over all drawings of $G$.
Clearly, (as non-crossing edges cross an even number of times)
for every graph $G$, 
$\ocr(G)\le\pcr(G)\le\cro(G)$.
According to the (weak) Hanani-Tutte theorem \cite{C34}, \cite{PSS07}, 
if $\ocr(G)=0$, then $G$ is planar, that is, $\ocr(G)=\pcr(G)=\cro(G)=0$.
It was shown in \cite{PSS07} that for $k=1, 2, 3$, if 
$\ocr(G)=k$, then $\ocr(G)=\pcr(G)=\cro(G)=k$.
There are examples where $\ocr$ is different from $\pcr$ and $\cro$,
namely, there is an infinite family of graphs with 
$\ocr(G)<0.855\cdot\pcr(G)$ \cite{T08}, \cite{PSS08}. From the other 
direction we only have $\pcr(G)<2\ocr(G)^2$ \cite{PT00a}.

In \cite{PT00b} some further variants were introduced, in order to study 
the role of crossings between adjacent edges. 
For each of $\cro$, $\pcr$, and $\ocr$, they introduced three counting rules.

\smallskip

\noindent {\bf Rule $+$:} Only those drawings are considered, where adjacent edges cannot cross.

\noindent {\bf Rule $0$:} Adjacent edges can  cross and their crossings are counted as well.

\noindent {\bf Rule $-$:} Adjacent edges can  cross and their crossings are not counted.

\smallskip

Combining these rules with the three crossing numbers, we get nine possibilities.
But 
it is easy to see that 
$\cro_{+}=\cro$ \cite{PT00b}.
On the other hand, regarding Rule $+$ for the odd-crossing number, 
it seems more natural to assume that adjacent edges cross an {\em even number of times} than
to assume that they do not cross at all. 
So, let $\ocr_{*}(G)$ be the minimum number of odd-crossing pairs of edges over all drawings of $G$ 
where adjacent edges cross an even number of times (these drawings are called {\em weakly semisimple} in \cite{BFK15}).
Therefore,  we have nine  versions, see Table 1.
In this table, values do not decrease if we move to the right or up, and it was shown in \cite{PSS08} that 
 $\cro(G)<2\ocr_{-}(G)^2$. On the other hand, there are graphs $G$ such that 
$\ocr_{-}(G)<\ocr(G)$ \cite{FPSS11}.

\medskip

\begin{center}
\begin{tabular}{|l|l|l|l|}
  \hline
  Rule $+$ & $\ocr_{*}(G)\le \ocr_{+}(G)$ & $\pcr_{+}(G)$ & \multirow{2}{*}{$\cro(G)$} \\
  \cline{1-3}
  Rule $0$ & $\ocr(G)$ & $\pcr(G)$ & \\
  \hline
  Rule $-$ & $\ocr_{-}(G)$ & $\pcr_{-}(G)$  & $\cro_{-}(G)$ \\
  \hline
\end{tabular}
\end{center}  

\smallskip














\centerline{{\bf Table 1.} Nine versions of the crossing number.}

\medskip

The Crossing Lemma, discovered by Ajtai, Chv\'atal, Newborn, Szemer\'edi 
\cite{ACNS82} and independently by Leighton \cite{L84} is 
definitely the most important inequality for crossing numbers.

\medskip

\noindent {\bf Crossing Lemma.} {\em If a simple graph $G$ of $n$ vertices has $m\ge 4.5n$ 
edges, then $\cro(G)\ge \frac{1}{60.75}\frac{m^3}{n^2}$ edges.}

\medskip

The bound is tight, apart from the value of the constant \cite{PT97}.
The constant above follows from the beautiful probabilistic argument of
Chazelle, Sharir and Welzl \cite{AZ04}. This argument works for all nine 
versions of the crossing number \cite{PT00b}.
For the classical crossing number, $\cro(G)$, the constant was improved 
in three steps \cite{PT97}, \cite{PRTT06}, 
the best bound is due to Ackerman \cite{A19},
$\cro(G)\ge \frac{1}{29}\frac{m^3}{n^2}$, when $m\ge 7n$.

The only improvement for any other version is 
a  result of Ackerman and Schaefer \cite{AS14},
$\pcr_{+}(G)\ge \frac{1}{34.2}\frac{m^3}{n^2}$, when $m\ge 6.75n$.
For all other versions of the crossing number, the constant $60.75$ is the
best we have.

In this note we get an improvement for two other versions, $\ocr_{+}$ and
$\ocr_{*}$.

\medskip

\begin{theorem}\label{oddcr}
  Suppose that $G$ has $n$ vertices and $m\ge 6n$ edges. Then
  $\ocr_{+}(G)\ge\ocr_{*}(G)\ge \frac{1}{54}\frac{m^3}{n^2}$.
\end{theorem}

Our approach is very similar to all previous improvements mentioned above.
The first step is to find many odd-crossing pairs in sparse graphs. 
Then this bound is applied for a random subgraph of $G$ to get the general bound. 







A graph $G$ is called {\em $k$-planar} if it can be drawn in the plane such that 
there are at most $k$ crossings on each edge. Such a drawing is called a {\em $k$-plane} drawing.
Let $m_k(n)$ denote the maximum number of edges of a $k$-planar graph of $n$
vertices.

Clearly, $m_0(n)=3n-6$.
It is known that  $m_1(n)=4n-8$ for $n\ge 12$,
$m_2(n)\le 5n-10$ and it is  tight for infinitely 
many values of $n$, \cite{PT97}, 
$m_3(n)\le 5.5n-11$, $m_4(n)\le 6n-12$, 
which are tight up to an additive constant
\cite{PRTT06}, \cite{A19}.

We prove an odd-even version of these results. 
A graph $G$ is called {\em $k$-odd-planar} if it can be drawn in the plane such that 
any edge is crossed {\em an odd number of times} by at most $k$ other edges
(edges crossing an even number of times do not count). Such a drawing is called a
 {\em $k$-odd-plane} drawing.

Let $m^{\mathrm{odd}}_k(n)$ denote the maximum number of edges of a
$k$-odd-planar 
graph with $n$ vertices.
Clearly, we have $m^{\mathrm{odd}}_k(n)\ge m_k(n)$
and by the weak Hanani-Tutte theorem 
 \cite{C34}, \cite{PSS07}, we have $m^{\mathrm{odd}}_0(n)=3n-6$.

\begin{theorem}\label{5n-9}
For any $n, k\ge 1$ we have
$$m^{\mathrm{odd}}_k(n)\le m_k(n)+k(n-1).$$
\end{theorem}




This result is interesing only for small $k$. For 
$k$ (and $n$) large enough, we have an easy better bound.

\begin{theorem}\label{gyokk}
For any $n, k\ge 1$ we have
$$m^{\mathrm{odd}}_k(n)\le \sqrt{32}\sqrt{k}n.$$
\end{theorem}

We do not think that our bounds are tight. 
We cannot even rule out the possibility that 
$m^{\mathrm{odd}}_k(n)=m_k(n)$ for every $n, k$.

\bigskip

\section{Proofs}

A (multi)graph $G$, together with its drawing $D$ in the plane, 
is called {\em topological (multi)graph}. 
The points (resp. curves), representing the vertices (resp. edges) 
of $G$ are called {\em vertices} (resp. {\em edges}) of the 
topological (multi)graph.
If the drawing is obvious from 
the context, we do not make any notational distinction 
between the topological (multi)graph and the underlying 
abstract (multi)graph. 
Let $G$ be a topological multigraph and $e$ an edge. 
The pieces of $e$ in small neighborhoods 
of its endpoints are called {\em endings} of $e$ and denoted by $e^+$ and $e^-$.
As $e$ is an undirected edge, the $+$ and $-$ signs have no special
meaning, either ending
can be  $e^+$ and the other one is $e^-$.
The {\em rotation system} 
is the cyclic order of adjacent edges, or endings, at each vertex.
A cyclic order is always clockwise.
Two edges form an {\em odd pair} 
(resp. {\em even pair}) if they cross an {\em odd} (resp. {\em even}) 
number of times.
An edge is called {\em even} if it is crossed
an even number of times by every other edge and it is {\em odd} otherwise.

\medskip

According to the weak Hanani-Tutte theorem, 
if a graph can be drawn so that any two edges cross an 
even number of times, then it is planar. This result has many proofs, 
one of the 
nicest and simplest is due to Pelsmajer, Schaefer 
and \v Stefankovi\v c \cite{PSS07}. The proof is based on the following lemma.

\smallskip

\noindent {\bf Lemma 0.  \cite{PSS07}} 
{\em Let $G$ be a topological multigraph that  has one vertex and $n$ edges (loops).
Suppose that every edge is even.
Then, $G$ can be redrawn such that the rotation system is the same and there is no 
edge crossing.}

\medskip

First we prove the following generalization. 

\medskip

\begin{lemma}\label{loop}
Let $G$ be a topological multigraph that  has one vertex and $m$ edges (loops).
Then, $G$ can be redrawn such that (i) the rotation system is the same
(ii) even pairs do not cross, (iii) odd pairs cross 
once, and (iv) there are no self-crossings.
\end{lemma}


\noindent {\it Proof.}
The proof is by induction on the number of edges.
If there is only one loop, the statement is trivial. 

Suppose that $G$ has one vertex $v$ and  
$m>1$ loops, 
and the statement 
has been proved for a smaller number of loops. Let $D$ denote the present drawing of $G$.
We can assume immediately that there are no self-crossings, since we can eliminate self-crossings of each loop
without changing the number of crossings between different loops. 

Let $s$ be a very short segment ending in $v$ and not crossing any of the
loops.  

Assume without loss of generality that
for any loop $e$, $s$, $e^-$ and $e^+$ are in this clockwise order at $v$
(otherwise we can switch the notation for $e^-$ and $e^+$).

For two loops $e$ and $f$, we say that $f\prec e$ if
the clockwise order of endings at $v$ is
$s$, $e^-$, $f^-$, $f^+$, $e^+$.
The relation $\prec$ defines a partial order on the loops.
Let $e$ be a minimal loop with respect to $\prec$.
That is, any other loop $f$ has at most one ending between $e^-$ and $e^+$. 

\begin{figure}
\begin{center}
\scalebox{0.35}{\includegraphics{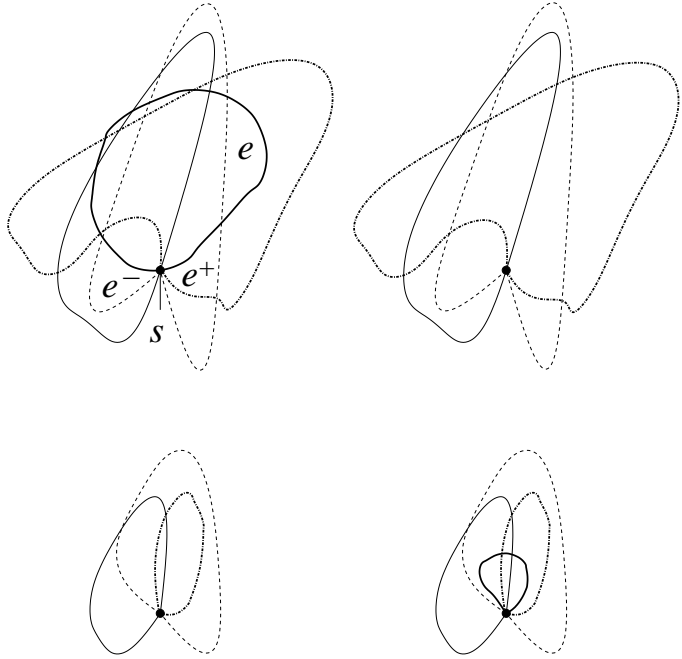}}
\caption{Steps in the proof of Lemma \ref{loop}.}\label{loopfigure}
\end{center} 
\end{figure}

Delete $e$ and 
apply the induction hypothesis.
We get a drawing $D'$ satisfying the conditions. In particular, the cyclic order of endings is the same as  
in $D$.
Insert $e^-$  and $e^+$ back to their original places and connect them by an arc in a small neighborhood of $v$, going clockwise
from $e^-$ to  $e^+$. We can do it such that it crosses only those endings, each exactly once,
which are between
$e^-$  and $e^+$.
In the obtained drawing $D_1$, the rotation system is the same as in $D$ and there are no self-crossings.
For any two loops  $f, g\neq e$, (ii) and (iii) are satisfied by the induction hypothesis.
If $e, f$ is an odd pair in $D$, then it has exactly one ending between $e^-$ and $e^+$ so in $D_1$
they cross once.
If it is an even pair in $D$, then $f$ has no 
ending between $e^-$ and $e^+$ so in $D_1$
they do not cross.
This concludes the proof of Lemma \ref{loop}.
See Figure \ref{loopfigure}.   \qed

\smallskip

\medskip

\noindent {\bf Remarks.} 1. The statement of Lemma \ref{loop} can be found implicitly in \cite{PSS07}, p. 492, as a remark.

\noindent 2. Another possible proof is the following. Take a  drawing of $G$ which has the same rotation system and under this condition
the minimum number of crossings. It can be shown that
this drawing satisfies the conditions,
otherwise we could have a better drawing of $G$.
But we did not find this method easier. 

\medskip

\begin{lemma}
Let $k\ge 0$, $l\ge 1$.
  Suppose that $n_1, \ldots, n_l>0$ and $n_1+\cdots +n_l=n$.
  Then

  $$m_k(n)\ge m_k(n_1)+\cdots +m_k(n_l).$$
\end{lemma}


\noindent {\em Proof.} For every $i$ take a $k$-plane graph of $n_i$ vertices and
$m_k(n_i)$
edges. Their disjoint union is a $k$-plane graph with $n$ vertices and 
$m_k(n_1)+\cdots +m_k(n_l)$
edges. \qed


\medskip

\begin{figure}[ht!]
\begin{center}
\scalebox{0.35}{\includegraphics{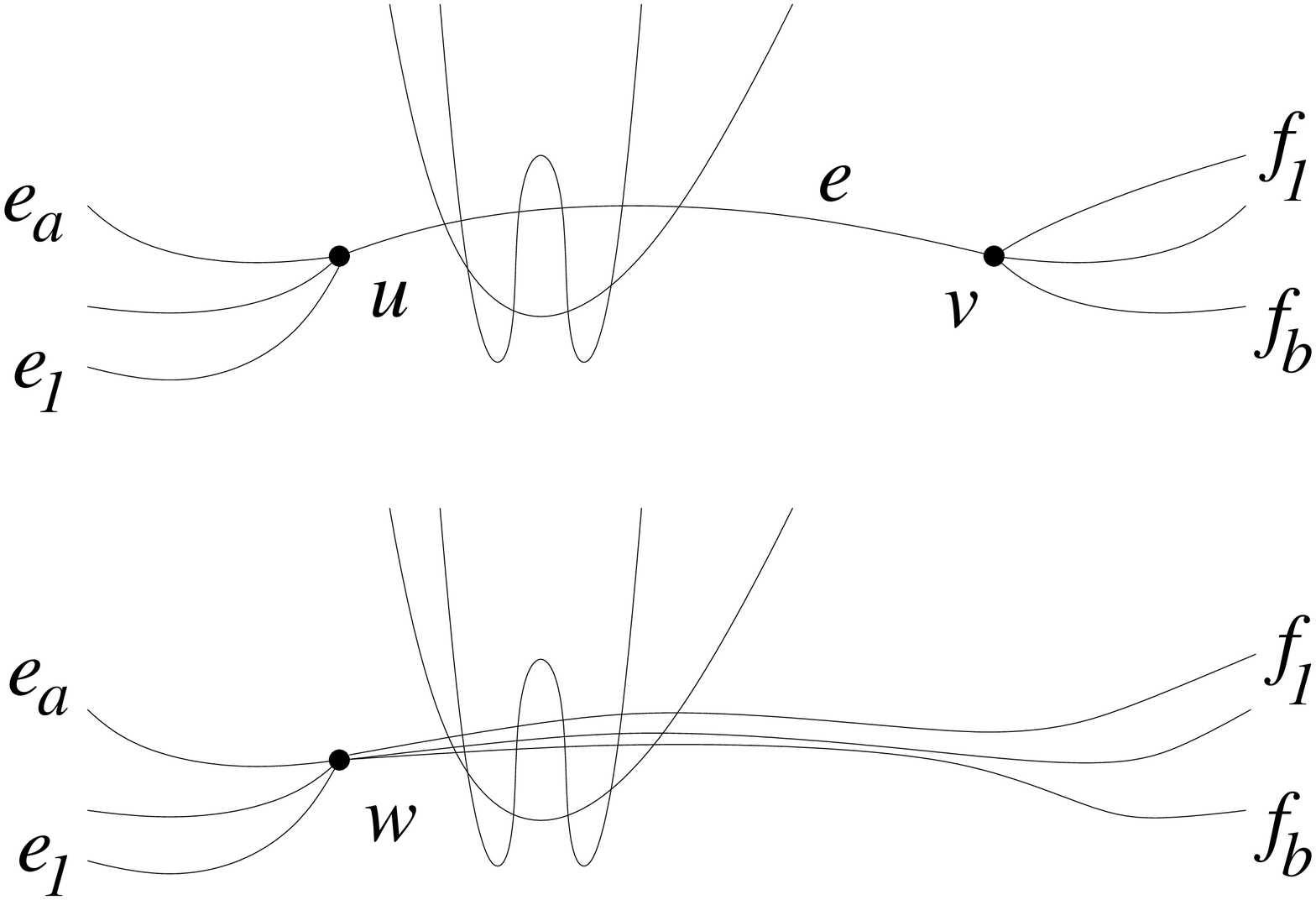}}
\caption{Contracting edge $e=uv$ to $u$. The resulting vertex is called $w$.}\label{contract}
\end{center}
\end{figure}

\noindent {\bf Definition.}
Suppose that $G$ is a topological multigraph. Let $e=uv$ be an edge.
The {\em contraction} of $e$ to $u$ is the following procedure. Move $v$ to $u$, along the edge $e$. 
At the same time, extend all other edges incident
to $v$ by curves along $e$, from a neighborhood of $v$ to $u$, without creating any crossing
among them. See Figure \ref{contract}.

\medskip

Suppose that $e=uv$ is an {\em even} edge of $G$. Contract $e$ to $u$ and let $w$ be the new vertex and $G'$ the resulting topological graph.
There is a natural bijection between the edges of $G\setminus e$ and $G'$, and the odd pairs are exactly the same.
Any vertex $x\neq u, v$ of $G$ remains a vertex of $G'$ with the same cyclic order of incident  edges.
Suppose that in $G$, $e, e_1, \ldots, e_a$ and $e, f_1, \ldots, f_b$ are the cyclic orders of incident
edges at $u$ and $v$, respectively.
Then in $G'$, the cyclic order of edges at $w$ is $e_1, \ldots, e_a, f_1, \ldots, f_b$.
We can ``reverse'' the contraction as follows.
If we replace $w$ in $G'$ by two very close  vertices $u'$ and $v'$, connected by an edge, $u'$ incident  to $e_1, \ldots, e_a$ and
$v'$ incident  to $f_1, \ldots, f_b$, we get a topological graph with exactly the same rotation system as $G$.

\medskip

\noindent {\bf Proof of Theorem \ref{5n-9}.} 
Suppose that $G$ has $n$ vertices and $m$ edges, and it 
is drawn in the plane such that any 
edge is crossed by at most $k$ other edges an odd number of times.
We can assume that $G$ is connected, as an abstract graph, otherwise, by Lemma 2, 
we can 
argue separately for each component.
Suppose first that $G$ contains a spanning tree $F$
whose edges cross each other an even number of times.
Remove all edges of $G$ that cross an edge 
of $F$ an odd number of times. We get the topological graph $G_1$ in the 
inherited drawing. 
Since $F$ has $n-1$ edges and each of them is crossed by  at most $k$ other 
edges an odd number of times, for the number of $m_1$ edges of $G_1$, $m_1\ge m-k(n-1)$.
In $G_1$, all edges of $F$ are even. Contract all edges of 
$F$. We get a topological multigraph $G_2$, which has only one vertex $v$. 
The odd pairs are exactly the same as in $G_1$.
Apply Lemma 1 for $G_2$ and get the topological multigraph $G_3$. Since the rotation system is the same in $G_3$ 
as in $G_2$, we can reverse the contractions without creating any additional crossing. 
This way we get $G_4$. Observe that $G_4$ is a redrawing of $G_1$,
with the property that if two edges form an even (odd) pair in $G$, then they do not cross (cross once) in $G_4$.
By the assumption, every edge in $G_1$ is part of at most $k$  odd pairs.
Therefore, 
$G_4$ is a $k$-plane drawing of $G_1$ so it has at most $m_k(n)$ edges.
Consequently, for the number of edges of $G$ we have $m\le m_k(n)+k(n-1)$.

Suppose now that $G$ does not 
contain a spanning tree 
whose edges cross each other an even number of times.
Let $F$ be a maximal forest whose edges cross each other an even number of times.
Let $V_1, \ldots, V_l$ be the vertices of the connected components of $F$, $|V_i|=n_i$.
Since $F$ is not a tree, $l\ge 2$ and $|E(F)|=n-l$.
By the maximality of $F$, those edges of $G$ that connect two components cross some edge of $F$ an odd number of times.
Remove all edges of $G$ that cross some edge of $F$ an odd number of times.
We removed at most $k(n-l)$ edges. The resulting graph has $l$ components $G_1, G_2, \ldots, G_l$, on vertices
 $V_1, \ldots, V_l$. Let $|E(G_i)|=m_i$. 
By the construction, each $G_i$ contains a spanning tree $F_i$ whose 
edges are even. Contract the edges of $F_i$ and argue as above. We get that
$$m\le \sum_{i=1}^lm_i+k(n-l)\le \sum_{i=1}^lm_k(n_i)+k(n-l)\le m_k(n)+k(n-1).$$
\qed


\medskip

\noindent {\bf Proof of Theorem \ref{oddcr}.} 
We have
$m^{\mathrm{odd}}_0(n)=3n-6$,
and by Theorem \ref{5n-9},
$m^{\mathrm{odd}}_1(n)\le m_1(n)+n-1=5n-9$.

First we show, 
by induction on the number of edges that
for any graph with $n$ vertices and $m$ edges,
$\ocr(G)\ge m-3n$. If $m\le 3n$, the statement is trivial.
Suppose that $m>3n$ and we have proved the statement for
$m-1$. Take any drawing $D$ of $G$. Since 
$m^{\mathrm{odd}}_0(n)=3n-6$, there is an odd pair $e, f$. Remove $e$ from $G$, and then
the obtained drawing has one less edges and at least one less odd pairs.
Therefore, the number of odd pairs in $D$ is at least $1+(m-1)-3n=m-3n$.

Now we show,  again by induction on the number of edges that
for any graph with $n$ vertices and $m$ edges,
$\ocr(G)\ge 2m-8n$. If $m\le 5n$, then $\ocr(G)\ge m-3n\ge 2m-8n$.
Suppose that $m>5n$ and we have proved the statement for
$m-1$. Take any drawing $D$ of $G$, since 
$m^{\mathrm{odd}}_1(n)\le 5n-9$, there is an edge $e$ in two odd pairs, $e, f$ and $e, g$.
Remove $e$ from $G$, and then  the obtained drawing has one less edges and at least two less odd pairs.
Therefore, the number of odd pairs in $D$ is at least $2+2(m-1)-8n=2m-8n$.


Let $G$ be a graph of $n$ vertices and $m\ge 6n$ edges,
drawn in the plane realizing $\ocr_{*}(G)$, that is, 
any two adjacent edges cross an  even number of times and there are 
$\ocr_{*}(G)$ pairs of edges that cross an odd number of times.
Take a random subgraph $G'$ such that we take each vertex independently with probability $p=6n/m$.
Let $n'$, $m'$, and $x(G')$ denote the number of vertices (resp. edges) of $G'$, and the 
number of odd-crossing pairs of edges in $G'$, in the inherited drawing.
We have 
$$E(n')=pn, \ \ E(m')=p^2m, \ \ E(\ocr_{*}(G'))\le E(x(G'))=p^4\ocr_{*}(G).$$
For $G'$ we have  $\ocr_{*}(G')\ge \ocr(G')\ge 2m'-8n'$, taking expected
values, 
$$p^4\ocr_{*}(G)\ge 2p^2m-8pn.$$
For  $p=6n/m$, we have 
$$\ocr_{*}(G)\ge\frac{1}{54}\frac{m^3}{n^2}.$$
\qed


\medskip

\noindent {\bf Remark.}
Combining Theorem \ref{5n-9} and the bounds for $m_k(n)$
we obtain that 
$m^{\mathrm{odd}}_1(n)\le 5n-9$ and 
$m^{\mathrm{odd}}_2(n)\le 7n-12$.
In the proof of
Theorem \ref{oddcr} we used only the first inequality, the second would not help.
However, if we could prove that $m^{\mathrm{odd}}_2(n)\le 6.8n+c$ for some constant $c$, 
then we would get an improvement in Theorem \ref{oddcr} as well.

\medskip

\noindent {\bf Proof of Theorem \ref{gyokk}.} 
We apply a version of the Crossing Lemma for the
odd-crossing number, from \cite{PT00a}:
If $G$ has $n$ vertices and $m$ edges, and $m\ge 4n$, then $\ocr(G)\ge \frac{1}{64}\frac{m^3}{n^2}$.

Let $k\ge 1$.  
Suppose that $G$ is $k$-odd-planar with 
$n$ vertices and  $m=m^{\mathrm{odd}}_k(n)$ edges.
We can assume that $m\ge 4n$, otherwise we are done. 

Take a $k$-odd-plane drawing of $G$.
Every edge participates in at most $k$ odd pairs. Therefore, there are
at most $km/2$ odd pairs, so $\ocr(G)\le km/2$.
On the other hand, we can apply the Crossing Lemma for the odd-crossing number, so
$$km/2\ge \ocr(G)\ge \frac{1}{64}\frac{m^3}{n^2}.$$
It follows that
$m\le  \sqrt{32}\sqrt{k}n$. \qed

\medskip

\noindent {\bf Remarks.}
1. It is easy to see that Theorem \ref{gyokk} is better than Theorem \ref{5n-9} if $n, k\ge 40$.
This threshold $40$ can be improved, but it is not very  interesting. 

\noindent 2. The analogue of Theorem \ref{gyokk} for $m_k(n)$ instead of $m^{\mathrm{odd}}_k(n)$ was first proved in \cite{PT97}.
Using the best known constant for the Crossing Lemma from \cite{A19} we can get that for any $n, k\ge 2$, $m_k(n)\le 3.81\sqrt{k}n$.

\noindent 3. In the proof of 
Theorem \ref{oddcr} we had to assume that adjacent edges
cross an even number of times. Therefore, it holds only for
$\ocr_{*}(G)$ and $\ocr_{+}(G)$. 
We believe that it is not necessary and 
Theorem \ref{oddcr} can be extended to 
$\ocr(G)$.
Similarly, we believe that the bound of Ackerman and Schaefer \cite{AS14},
$\pcr_{+}(G)\ge \frac{1}{34.2}\frac{m^3}{n^2}$ can be extended to $\pcr(G)$.
We also believe that both of them are very interesing problems.

\end{document}